\documentclass[12pt]{amsart}

\input pictex
\usepackage[mathscr]{euscript}
\theoremstyle{plain}
\newtheorem{theorem}{Theorem}[section]
\newtheorem{lemma}[theorem]{Lemma}

\theoremstyle{definition}
\newtheorem{example}[theorem]{Example}
\newtheorem{remark}[theorem]{Remark}
\numberwithin{equation}{section}

\newcommand{\C}{\mathbb C}
\newcommand{\R}{\mathbb R}
\renewcommand{\P}{\mathbb P}

\newcommand{\spn}{\operatorname{span}}
\newcommand{\pa}{\partial}
\newcommand{\zb}{\overline z}
\newcommand{\Wp}{W^\prime}
\newcommand{\Sh}{\widehat S}
\begin{document}

\title[Hypersurface Edge-of-Wedge]
{An Edge-of-the-Wedge Theorem for\\[8pt] Hypersurface CR Functions}
\author[M.G. Eastwood]{Michael G. Eastwood}
\address{
  Department of Pure Mathematics\\
  University of Adelaide\\
  South AUSTRALIA 5005}
\email{meastwoo@maths.adelaide.edu.au}

\author[C.R. Graham]{C. Robin Graham}
\address{
  Department of Mathematics\\
  University of Washington \\
  Box 354350\\
  Seattle, WA 98195-4350}
\email{robin@math.washington.edu}
\maketitle

\thispagestyle{empty}

\renewcommand{\thefootnote}{}
\footnotetext{This research was supported by the Australian Research Council
and the University of Washington. This support and the hospitality of the
Universities of Adelaide and Washington is gratefully acknowledged.}

\section{Introduction}\label{intro}
The Lewy extension theorem asserts the holomorphic extendability of CR
functions defined in a neighborhood of a point on a hypersurface in
$\C^{n+1}$.  The edge-of-the-wedge theorem asserts the extendability of
holomorphic functions defined in wedges in $\C^{n+1}$ with edge a maximally
real submanifold.  In this article we prove under suitable hypotheses the
holomorphic extendability to an open set in $\C^{n+1}$ of CR functions
defined in the intersection of a hypersurface with a wedge whose edge is
contained in the hypersurface.  Unlike the situation for the classical
edge-of-the-wedge theorem, for this hypersurface version extendability
generally depends on the direction of the wedge.

Equip $\R^{n+1}$ with its standard inner product and let
${\sigma}\in\R^{n+1}$ be
a unit vector.  By the {\em round cone\/} in $\R^{n+1}$ of aperture
$\delta>0$, extent $\ell>0$, and axis ${\sigma}$, we shall mean
$$\{x\in\R^{n+1}:|x-\langle x,{\sigma}\rangle
{\sigma}|<\delta\langle x,{\sigma}\rangle \,\mbox{ and }|x|<\ell\}.$$
Let $E\subset\C^{n+1}$ be a {\em maximally real\/}
submanifold, i.e. $E$ is a totally real submanifold of maximal dimension
$n+1$.  Using $J$ to denote multiplication by~$i$ and $TE$ to denote the
tangent bundle to $E$, for $p \in E$ the subspace  $J(T_pE)$ is transverse
to~$T_pE$.  Now $J(T_pE)$ inherits an inner product from that on
$\C^{n+1}$, hence may be identified
with $\R^{n+1}$ so that round cones in $J(T_pE)$ are defined.  We may view
$J(T_pE)$ as a real affine subspace of $\C^{n+1}$.  Shrinking
$E$ if necessary, we assume that for some $r > 0$ the
balls of radius $r$ in the subspaces $J(T_pE)$ sweep out a
tubular neighborhood of $E$ in $\C^{n+1}$.  By a {\em one-sided wedge\/}
in $\C^{n+1}$
of aperture $\delta>0$, extent $\ell\in(0,r)$, and  edge $E$
we shall mean the union of the round cones of aperture~$\delta$, extent
$\ell$, and axes ${\sigma}(p)$ for some smooth family of unit
vectors ${\sigma}(p) \in J(T_pE)$ for~$p\in E$.  By the opposite of such a
one-sided wedge in $\C^{n+1}$  we shall mean
the union of the round cones of aperture~$\delta$, extent
$\ell$, and axes $-{\sigma}(p)$, and by a {\em two-sided wedge\/} in
$\C^{n+1}$
we shall mean the union of such a one-sided wedge with its opposite.
The classical edge-of-the-wedge
theorem (e.g.~\cite{r}) asserts that a holomorphic function defined on a
two-sided wedge in $\C^{n+1}$ and continuous across its edge extends as a
holomorphic function to a neighborhood of this edge.

Now suppose that $M$ is a smooth hypersurface in $\C^{n+1}$.
The holomorphic tangent bundle $H \equiv TM\cap J(TM)$ on $M$
together with $J|_H$ defines the CR structure of $M$ intrinsically.
Let $E\subset M$ be a
smooth submanifold which is maximally real as a submanifold of $\C^{n+1}$.
Define a subbundle $N \subset TM|_E$ by $N = J(TE) \cap TM$.  Observe that
$N \subset H$ and for $p \in E$, $N_p$ is a complement to $T_pE$ in
$T_pM$.  By a one-sided wedge in~$M$ with edge $E$ we shall mean
the intersection $W^+$ of $M$ with a one-sided wedge in $\C^{n+1}$ with
edge $E$ and with axes ${\sigma}(p)$ everywhere tangent to $M$.
We denote by
$W^-$ the  intersection of $M$ with the opposite one-sided wedge in
$\C^{n+1}$.  By a two-sided wedge in $M$ with edge $E$ we shall mean the
union $W=W^+ \cup W^-$ of two such opposite one-sided wedges in $M$.
Observe that the axes of a
wedge in $M$ satisfy ${\sigma}(p) \in N_p$.  We shall say that a
non-zero vector
$\tau \in N_p$ is a {\em direction\/} of $W^+$ (resp.\ $W^-$) at $p$
if $\tau$ is in the intersection of $T_pM$
with the round cone of infinite extent in $J(T_pE)$ defining $W^+$
(resp.\ $W^-$).  We say that ${\tau}$ is a direction of $W$ if it is a
direction of either $W^+$ or~$W^-$.

Locally we may choose a smooth real defining function $r$ for $M$ so that
$M=\{r=0\}$ and $dr \neq 0$ everywhere.  The
{\em Levi form\/} associated to
$r$ is the complex-valued quadratic form on $H$,
Hermitian with respect to $J$, given by
$L({\sigma},{\tau})=\partial \overline{\partial} r
( {\sigma}-iJ{\sigma}, {\tau}+iJ{\tau} )$ for ${\sigma},{\tau} \in H$.
For $p \in M$ we say that $L_p$ is {\em indefinite\/} if there are
${\sigma}_+, {\sigma}_- \in H_p$ such that $L({\sigma}_+,{\sigma}_+) > 0$
and $L({\sigma}_-,{\sigma}_-) <0$.  We say
that ${\sigma} \in H_p$ is {\em null\/} if $L({\sigma},{\sigma}) = 0$.

Our main theorem is as follows.

\begin{theorem}\label{main}
Let $M$ be a smooth hypersurface in
$\C^{n+1}$ and $E\subset M$ a maximally real submanifold. Let $W$ be a
two-sided wedge
in $M$ with edge~$E$. Suppose that $p_0 \in E$ is a point at which the Levi
form
of $M$ is indefinite and that $W$ has a direction ${\sigma} \in N_{p_0}$
which is null. Then there is a neighborhood of $p_0$ in $\C^{n+1}$ to which
every CR function on $W$,
continuous across $E$, extends as a holomorphic function.
\end{theorem}

Note that by replacing the round cones defining $W^+$ by smaller ones, we
can suppose, without loss of generality, that the axis ${\sigma}(p_0)$
is null.  Note also that the Levi form $L_p$ is indefinite for all
$p \in M$ near $p_0$.  Therefore by the Lewy extension theorem applied
to points $p \neq p_0$, any CR
function on $W$ is, close enough to $p_0$, the restriction of a holomorphic
function, so in particular is a smooth function.  The Lewy extension
theorem also allows us to rephrase the conclusion of Theorem~\ref{main}
intrinsically by saying that there is a neighborhood of $p_0$ in $M$ to
which we obtain extension as a CR function.

As indicated earlier, Theorem~\ref{main} generalizes at the same time the
(two-sided) Lewy extension theorem and the edge-of-the-wedge theorem.  In
the Lewy theorem, the CR function is required to be defined in a full
neighborhood in $M$ of $p_0$ in order to get extension; in the
edge-of-the-wedge theorem, the holomorphic function is required to be
defined in a full two-sided wedge in $\C^{n+1}$.

Since Theorem~\ref{main} requires the Levi form to be indefinite,
the hypotheses cannot be satisfied unless $n \geq 2$.  The hypersurface
edge-of-the-wedge phenomena of interest to us do not occur in $\C^2$ in any
case since a maximally real submanifold $E$ of a hypersurface $M$ of
$\C^2$ has codimension 1, so a two-sided wedge in $M$ together with its
edge $E$ actually fills out a full neighborhood of $E$.  We consider two
illustrative examples in $\C^3$.
Let $(z_1,z_2,w) =(x_1 + iy_1,x_2+iy_2,u+iv)$
denote the co\"ordinate functions on $\C^3$.

\begin{example}\label{counterexample}
Let $M=\{v=y_1^2-y_2^2\}$ and $E=\R^3=\{y=v=0\}$.
Then $H=\spn\{\pa_{x_1},\pa_{x_2},\pa_{y_1}, \pa_{y_2}\}$ and
$N=\spn\{\pa_{y_1},\pa_{y_2}\}$.
The Levi-null directions in $N$ are multiples of
$\pa_{y_1} \pm \pa_{y_2}$,
so  Theorem~\ref{main} implies that if one of these is a
direction of $W$, then extension must hold.  On the other hand,
any non-null direction in $N$ is a direction of a two-sided wedge $W$ in
$M$ for which there do exist CR functions on $W^{\pm}$
with equal boundary values on $E$ for which extension fails.  To see this,
consider $W=M\cap\{v>0\}=W^+\cup W^-$, where
$W^{\pm}=M\cap\{\pm y_1>|y_2|\}$.  Sufficiently near $E$, $W^+$ is the
intersection of $M$ with the union over points of $E$ of the round cones of
aperture 1 and axes in the $y_1$-direction.  Consequently,
the intersection of $W$ with a small enough neighborhood of $E$ is a
two-sided wedge in $M$ with edge $E$, having all directions
$a \pa_{y_1} + b \pa_{y_2}$ with $|a| > |b|$.
Let $f$ be a holomorphic function in the upper half plane in $\C$
which extends smoothly to the closed upper half plane but which does not
have a holomorphic extension across $0$.  Then
$f(w)|_W$ defines CR functions on $W^{\pm}$ with
the required properties.  Clearly a similar construction can be carried out
for $M\cap\{v<0\}$.
\end{example}

\begin{example}\label{hyperquadric}
Let $M=\{v=\Im (z_1 \zb_2)=y_1x_2-y_2x_1\}$, and again let
$E=\R^3$.  This $M$ is
the standard mixed signature hyperquadric in $\C\P_3$ viewed in an affine
chart.  In this case $L$ is indefinite
and every direction $\sigma \in N$ is null, so extension holds for any
wedge.
We obtain a result for CR functions completely analogous to the classical
edge-of-the-wedge theorem for holomorphic functions.
\end{example}

For the case of two-sided wedges in $M$ which point in directions which are
not null, we have the following one-sided extension theorem.

\begin{theorem}\label{one-sided}
Let $M$ be a smooth hypersurface in $\C^{n+1}$ with defining function $r$
and let  $E\subset M$ be a maximally real submanifold. Let $W$ be a
two-sided wedge in $M$ with edge~$E$.  Suppose that $p_0 \in E$ and that
$W$ has a direction $\sigma \in N_{p_0}$ satisfying $L(\sigma,\sigma)>0$.
Then there is an open set $U\subset \{r<0\}$ with the following properties:
\begin{enumerate}
\item\label{firstitem}
For each $p \in E$ sufficiently close to $p_0$ and each unit vector $\tau \in
J(T_pE)$ such that $dr(\tau)<0$, there is a round cone in $J(T_pE)$ with axis
$\tau$ which is contained in $U$.

\item\label{seconditem}
$\overline U$ contains a two-sided subwedge $\widetilde W$ in $M$ of $W$ having
$\sigma$ as a direction at~$p_0$.

\item To every continuous CR function $f$ on $W$, continuous across $E$, there
is a holomorphic function on $U$, continuous in $\overline U$, which restricts
to $f$ on $\widetilde W$.  In particular, the extension agrees with $f$ on~$E$.
\end{enumerate}
\end{theorem}

We remark that
Theorem~\ref{one-sided} can be applied to a two-sided wedge in $M$ with
direction $\sigma$ satisfying $L(\sigma,\sigma)<0$ upon replacing $r$ by
$-r$.  We also remark that in \ref{firstitem} above, the aperture and extent
of the cone with axis $\tau$ may be chosen to depend continuously on $p$ and
$\tau$, but they will in general go to $0$ as $dr(\tau) \rightarrow 0$.

Theorem~\ref{one-sided} may be seen as a refinement of one-sided
Lewy extension in which the CR function need only be defined in a wedge in
$M$.  Example~\ref{counterexample} shows that extension to all of $\{r<0\}$
need not hold in this case.

In this paper we prove Theorems~\ref{main} and \ref{one-sided} by analytic
disc methods combined with an adaptation of the Baouendi-Treves
approximation theorem to wedges.  Such a proof of the
$\C^{n+1}$-edge-of-the-wedge  theorem for $C^1$ edges was given by
Rosay~\cite{r}.
The discs required in Theorem~\ref{one-sided} are direct
modifications of the usual discs used in Lewy extension.  For
Theorem~\ref{main} we first use Lewy extension away from the edge to obtain
an open set in $\C^{n+1}$ containing $W$; this set has a parabolic
``spike'' approach transverse to $M$.  We then use a version of the
``folding screen lemma'' to enlarge each side of the wedge to a full wedge
in $\C^{n+1}$, whereupon we can apply the usual edge-of-the-wedge theorem.
We remark that if in Theorem~\ref{main}, $L|_{N_{p_0}}$ is
indefinite in the sense that there are
$\tau_+, \tau_- \in N_{p_0}$ with
$L({\tau}_+,{\tau}_+) > 0$ and $L({\tau}_-,{\tau}_-) <0$, then a proof can
be given using Theorem~\ref{one-sided}.  In fact, in this case any null
$\sigma \in N_{p_0}$ can be perturbed to $\sigma_+, \sigma_- \in
N_{p_0}$ satisfying $L({\sigma}_+,{\sigma}_+) > 0$ and
$L({\sigma}_-,{\sigma}_-) <0$, so Theorem~\ref{one-sided} gives
extension to one-sided wedges in $\C^{n+1}$ on either side of $M$.  The
usual edge-of-the-wedge theorem then yields holomorphic extension to a
neighborhood of $p_0$.  This argument can be used for
Example~\ref{counterexample}.  However, it seems that
Example~\ref{hyperquadric} requires a more sophisticated proof.

In \cite{t}, Tumanov considers extension of a CR function from a
one-sided wedge. In our situation his results imply that there is a
holomorphic extension to an ambient one-sided wedge but give no information
on the directions or axis of this ambient wedge. His results apply equally
to Example~\ref{hyperquadric}, for which the intrinsic edge-of-the-wedge
theorem holds, and to Example~\ref{counterexample}, for which it generally
fails.  In order to apply the ambient classical
edge-of-the-wedge theorem as we do in Theorem~\ref{main}, 
it is crucial that we obtain holomorphic
extension to a wedge with the same axis as the intrinsic wedge (as detailed in
Remark~\ref{sameaxis}).

In another paper \cite{micro} we will present a microlocal approach to
results of this type formulated in terms of the hypo-analytic wave-front
set of~\cite{bct}.  This in particular will contain analogues of
Theorems~\ref{main} and \ref{one-sided} for higher codimension CR
manifolds, and, in fact, an intrinsic version on general hypo-analytic
manifolds.

Throughout this article, smooth will mean infinitely differentiable.

\section{Folding Screens}\label{spikes}
In \cite{h}, H\"ormander
proved Bochner's tube theorem by a geometrical
arrangement dubbed a `folding screen' by Komatsu~\cite{k}. We shall use this
arrangement to study the polynomial hull of certain sets of wiggling spikes
in~$\C^2$. If $S \subset \C^{n+1}$, we denote by $\Sh$ the polynomial hull
of $S$: $$\widehat S = \{z: |p(z)| \leq \sup_{\zeta \in S}|p(\zeta)|
\,\,\mbox{for all holomorphic polynomials}\,\, p \}.$$

\begin{lemma}\label{preliminary_folding_screen}
Set
$$S=\left\{\begin{array}l
                (\zeta_1,\zeta_2)=(\xi_1+i\eta_1,\xi_2+i\eta_2)\in\C^2:\\[3pt]
\phantom{(\zeta_1,\zeta_2)=} 0<\eta_1<2
\mbox{\rm\ and }|\eta_2|<2\eta_1{}^2-|\xi|^2\eta_1/4
\end{array}\right\}.$$
Then $\Sh$ contains
$$T\equiv\left\{(\zeta_1,\zeta_2)=(i\eta_1,i\eta_2):
0<\eta_1<1\mbox{\rm\ and }|\eta_2|<\eta_1/2\right\}.$$
\end{lemma}
\begin{proof}
We construct analytic discs with boundaries contained in $S$.
For each $t\in(5,6)$, we claim
that the non-singular quadric
$$\{(\zeta_1+i)^2+(\zeta_2-2i)^2+t=0\}$$
meets the square pillar
$$\{(\zeta_1,\zeta_2)=(\xi_1+i\eta_1,\xi_2+i\eta_2):
                             0<\eta_1,\eta_2<1\}$$
in an analytic disc $D_t$
with piecewise smooth boundary in the folding screen
$$\{\xi + i\eta : 0<\eta_1\leq 1,\, \eta_2=0\}\cup
                     \{\xi + i\eta : \eta_1=1,\, 0\leq\eta_2<1\}.$$
In fact, the disc in question, with its boundary, may be
parameterized by
$$\{z=x+iy: t-4+y^2/16 \leq x \leq 4-y^2/16 \} \subset \C,$$
for it is easy to check that any such $z$ may be written uniquely both as
$$-(\zeta-2i)^2\mbox{ for some }\zeta=\xi+i\eta\mbox{ with }0\leq\eta<1$$
and
$$(\zeta+i)^2+t\mbox{ for some }\zeta=\xi+i\eta\mbox{ with }0<\eta\leq 1.$$
Since
$$\Re\left[(\zeta_1+i)^2+(\zeta_2-2i)^2+t\right]
=|\xi|^2+t-\left[(\eta_1+1)^2+(\eta_2-2)^2\right],$$
it follows that the disc $D_t$ is contained in the hypersurface
\begin{equation}\label{disc_hypersurface}
\{(\eta_1+1)^2+(\eta_2-2)^2=|\xi|^2+t\}\end{equation}
which intersects the folding screen as
$$\begin{array}l\{0<\eta_1=\sqrt{|\xi|^2+t-4}\; -1\leq 1,\, \eta_2=0\}\;
\cup\\[3pt]
\hspace*{100pt}
\{\eta_1=1,\, 0\leq\eta_2=2-\sqrt{|\xi|^2+t-4}<1\}.\end{array}$$
This intersection is empty unless $|\xi|^2\leq 3$ and, in this case, it is
easily verified that
$$\begin{array}l\phantom{\mbox{and\quad}}
|\xi|^2/8\leq\sqrt{|\xi|^2+1}\; -1<\sqrt{|\xi|^2+t-4}\; -1\\[3pt]
\mbox{and\quad}2-|\xi|^2/4>2-\sqrt{|\xi|^2+1}>2-\sqrt{|\xi|^2+t-4}.
\end{array}$$
Since $\partial D_t$ is the intersection of $\overline{D_t}$
with the folding screen,
from these inequalities it follows that $\partial D_t\subset S$ for all
$t\in(5,6)$.

The following diagram shows
what is happening in the $(\eta_1,\eta_2)$-plane when~$\xi=0$.
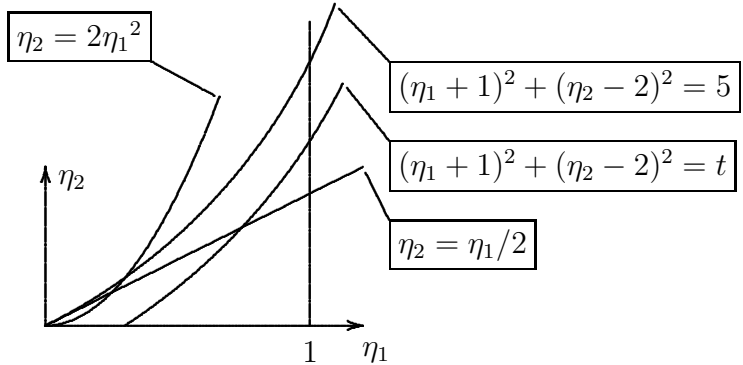
\begin{center}
\mbox{\beginpicture
\setplotarea x from 0 to 100, y from -20 to 125
\arrow <6pt> [.2,.6] from 0 0 to 120 0
\arrow <6pt> [.2,.6] from 0 0 to 0 60
\put {$\eta_1$} at 125 -10
\put {$\eta_2$} at 10 55
\put {$1$} at 100 -10
\plot 0 0  120 60 /
\plot 100 0  100 115 /
\plot 130 100  111 120 /
\put {\framebox{$(\eta_1+1)^2+(\eta_2-2)^2=5$}} [lt] at 130 100
\plot 130 70  114 90 /
\put {\framebox{$(\eta_1+1)^2+(\eta_2-2)^2=t$}} [lt] at 130 70
\plot 130 40  121 58 /
\put {\framebox{$\eta_2=\eta_1/2$}} [lt] at 130 40
\plot 40 100  64 87 /
\put {\framebox{$\eta_2=2\eta_1{}^2$}} [rb] at 40 100
\setquadratic
\plot 0 0  11 2.42  22 9.68  33 21.78  44 38.72  55 60.50  66 87.12 /
\circulararc 43 degrees from 0 0 center at -100 200
\circulararc 30 degrees from 30 0 center at -100 200
\endpicture}\end{center}
It is clear that every point in $T \cap \{\eta_2>0\}$,
not already in $S$, lies in $D_t$ for some $t\in(5,6)$.
A similar family of discs sweeps out $T \cap \{\eta_2<0\}$.
The result thus follows from the maximum modulus principle.
\end{proof}

\begin{lemma}\label{folding_screen}
For any $\epsilon,K>0$, there is a $\delta>0$ such that if we set
$$S_{\epsilon,K}=\left\{\begin{array}l
                (\zeta_1,\zeta_2)=(\xi_1+i\eta_1,\xi_2+i\eta_2)\in\C^2:\\[3pt]
\phantom{(\zeta_1,\zeta_2)=} 0<\eta_1<\epsilon
\mbox{\rm\ and }|\eta_2|<\epsilon\eta_1{}^2-K|\xi|^2\eta_1
\end{array}\right\},$$
then $\widehat{S_{\epsilon,K}}$ contains

$$T_{\delta}\equiv\left\{(\zeta_1,\zeta_2)=(i\eta_1,i\eta_2):
0<\eta_1<\delta\mbox{\rm\ and }|\eta_2|<\delta\eta_1\right\}.$$
\end{lemma}
\begin{proof}
Consider the linear change of co\"ordinates:
$$\phi(\zeta_1,\zeta_2) = (A\zeta_1, 2A^2\zeta_2/\epsilon) =
(\hat\zeta_1,\hat\zeta_2),$$
where $A=\max\{2/\epsilon,\epsilon/2,8K/\epsilon\}.$
We claim that $\phi^{-1}(S) \subset ~S_{\epsilon,K}$.
Were this to be shown, the result would follow
immediately from Lemma~\ref{preliminary_folding_screen}. Since
$A\geq 2/\epsilon$,
$$0<\hat\eta_1<2 \Rightarrow 0<\eta_1<\epsilon.$$
Since $A\geq\epsilon/2$,
$$|\hat\xi|^2=(A\xi_1)^2+(2A^2\xi_2/\epsilon)^2\geq
(A\xi_1)^2+(A\xi_2)^2=A^2|\xi|^2.$$
Thus,
$$\begin{array}{rcl}|\hat\eta_2|<2\hat\eta_1{}^2-|\hat\xi|^2\hat\eta_1/4
&\Rightarrow&
|2A^2\eta_2/\epsilon|<2(A\eta_1)^2-A^2|\xi|^2(A\eta_1)/4\\[3pt]
&\Rightarrow&
|\eta_2|<\epsilon\eta_1{}^2-K|\xi|^2\eta_1
\end{array}$$
since $A\geq 8K/\epsilon$. Altogether, $\phi^{-1}(S)\subset S_{\epsilon,K}$
as required. \end{proof}

We define the spike $Sp(\beta,\ell,0,0,m) \subset \R^2$ of sharpness $\beta
>0$, length $\ell >0$, slope $m \in \R$, and with vertex at the origin, by
$$Sp(\beta,\ell,0,0,m)\equiv
\left\{
\begin{array}l(\eta_1,\eta_2)\in\R^2 : 0<\eta_1, \\[3pt]
\phantom{(\eta_1,\eta_2)}
|(\eta_1,\eta_2)|<\ell, \mbox{ and }
|\eta_2-m \eta_1|<\beta \eta_1{}^2
\end{array}\right\}.$$
Let $Sp(\beta,\ell,q_1,q_2,m)$ denote the translation of this spike with
vertex at $(q_1,q_2) \in \R^2$.

\begin{lemma}\label{spike_lemma}
Suppose that $A,\beta,\ell,r$ are positive constants and suppose
$$S\supset
\left\{\begin{array}l(\zeta_1,\zeta_2)=(\xi_1+i\eta_1,\xi_2+i\eta_2):
|\xi|<r\;\mbox{\rm\ and}\\[3pt]
\phantom{(\zeta_1,\zeta_2)={}}(\eta_1,\eta_2)\in
Sp(\beta,\ell,q_1(\xi),q_2(\xi),m(\xi))
\end{array}\right\},$$
where $q_1(\xi),q_2(\xi),m(\xi)$ are smooth functions satisfying
$$|q_1(\xi)|\leq A|\xi|^2\qquad|q_2(\xi)|\leq A|\xi|^4
\qquad|m(\xi)|\leq A|\xi|^2.$$
Then there exists $\delta>0$ depending only on $A,\beta,\ell,r$ such that
$T_{\delta} \subset \Sh$.
\end{lemma}
\begin{proof}
The result follows from Lemma~\ref{folding_screen} provided we can find
$\epsilon,K>0$ so that $S \supset S_{\epsilon,K}$. Suppose
$\zeta=(\zeta_1,\zeta_2)=(\xi_1+i\eta_1,\xi_2+i\eta_2) \in S_{\epsilon,K}$.
Then $0<\eta_1<\epsilon$ and $\epsilon\eta_1{}^2-K|\xi|^2\eta_1>0$, so
$|\xi|<\epsilon/\sqrt{K}$. If we take $\epsilon\leq 1$ and $K\geq 1/r^2$,
then $|\xi|<r$ which is one of the conditions forcing $\zeta$ to be
in~$S$.  If we also take $K\geq A$, then
$$\eta_1>K|\xi|^2/\epsilon\geq A|\xi|^2\geq q_1(\xi),$$
which is the first condition for
$(\eta_1,\eta_2)\in Sp(\beta,\ell,q_1(\xi),q_2(\xi),m(\xi))$.
Now $|\eta_2|<\epsilon\eta_1{}^2$. Therefore,
$|\eta|\leq2\epsilon$. Thus,
$$\begin{array}{rcl}
|\eta-(q_1(\xi),q_2(\xi))|&\leq&|\eta|+|q_1(\xi)|+|q_2(\xi)|\\[3pt]
&\leq&2\epsilon+A|\xi|^2+A|\xi|^4\\[3pt]
&\leq&2\epsilon+A(1/K+1/K^2)
\end{array}$$
which we can arrange to be less than $\ell$ by taking $\epsilon$
sufficiently small and $K$ sufficiently large.
This leaves one condition to be ensured.  It is that
$$\beta(\eta_1-q_1(\xi))^2-|\eta_2-q_2(\xi)-m(\xi)(\eta_1-q_1(\xi))|>0$$
when $\zeta=\xi+i\eta\in S_{\epsilon,K}$. Now
$|\eta_2|<\epsilon\eta_1{}^2-K|\xi|^2\eta_1$ on $S_{\epsilon,K}$ so
it suffices to show that
$$\beta(\eta_1-q_1(\xi))^2-\epsilon\eta_1{}^2+K|\xi|^2\eta_1
-|q_2(\xi)+m(\xi)(\eta_1-q_1(\xi))|>0
$$
when $\eta_1>K|\xi|^2/\epsilon$. Choosing $\epsilon<\beta/2$, this
expression is strictly bounded below by
$$\beta\eta_1{}^2/2+(K|\xi|^2-2\beta q_1(\xi)-|m(\xi)|)\eta_1
-|q_2(\xi)|-|m(\xi)q_1(\xi)|$$
which, for $\eta_1>K|\xi|^2/\epsilon$, is bounded below by
$$\beta K^2|\xi|^4/2\epsilon^2+
(K|\xi|^2-2\beta A|\xi|^2-A|\xi|^2)\eta_1
-A|\xi|^4-A^2|\xi|^4.$$
Taking $K>(2\beta+1)A$, we may neglect the second term, leaving
$$(\beta K^2/2\epsilon^2
-A-A^2)|\xi|^4,$$
which is non-negative if we choose $K^2/\epsilon^2\geq 2(A+A^2)/\beta$.
\end{proof}

\section{A Normal Form}
The following normal form was suggested to us by Vladimir Ezhov.
\begin{lemma}\label{normal}
Suppose that $M$ is a smooth real hypersurface in~$\C^{n+1}$,
$E \subset M$ is a maximally real submanifold, and $p \in E$.
Then we may choose holomorphic co\"ordinates near $p$
$$(z_1,z_2,\dots,z_n,w)=(x_1+iy_1,x_2+iy_2,\dots,x_n+iy_n,u+iv)$$
so that $z=w=0$ at $p$,
$M$ has a defining function with Taylor series
$$r=-v+y^t\Lambda y+x^t\Omega y
         +\mbox{\rm cubic and higher order terms in $x$, $y$, and $u$},$$
and $E$ osculates $\R^{n+1}$ to order three at the origin, i.e.~$E$ may be
written near the origin as the graph of a smooth function
$\R^{n+1}\to i\R^{n+1}$ whose Taylor series begins with terms of order at
least four.
Here, $\Lambda $ is a symmetric $n\times n$ matrix and $\Omega$ is a skew
$n\times n$ matrix.
\end{lemma}
\begin{proof}
We first translate $p$ to the origin, and then by a complex linear
transformation, we may assume that
\begin{equation}\label{tangent_spaces}
T_0M=\{v=0\}\quad\mbox{and}\quad T_0E=\{y_1=y_2=\cdots=y_n=v=0\}.
\end{equation}
Then, near the origin,
$$E=\{y_1=f_1(x,u),y_2=f_2(x,u),\dots,y_n=f_n(x,u),v=g(x,u)\}$$
for smooth functions $f(x,u)$ and $g(x,u)$ whose Taylor series begin with
terms of order at least two.
Let $F(x,u)$ and $G(x,u)$ denote their third order
Taylor polynomials. Then the complex polynomial change of co\"ordinates
$$z=\hat z+iF(\hat z,\hat w)\qquad w=\hat w+iG(\hat z,\hat w)$$
is invertible near the origin and, in these new co\"ordinates, $E$ osculates
$\R^{n+1}$ to third order, as required. We shall now assume that this is done
and drop the hats. The equation for $M$ takes the form
$$2v=z^t(\Lambda +i\Omega)\bar z-\Re(z^t\Lambda z)
+\Im(z^t{\Gamma}z)+uy^t\mu+\cdots$$
where $L=\Lambda +i\Omega$ is an $n\times n$ Hermitian matrix (the Levi
form), ${\Gamma}$ is
a real symmetric matrix, $\mu$ is a real vector, and the ellipsis $\cdots$
indicates cubic and higher order terms. The change of co\"ordinates
$$w\mapsto w+\textstyle\frac{1}{2}(z^t{\Gamma}z+wz^t\mu)$$
preserves $\R^{n+1}$ and gives a new equation
$$2v=z^t(\Lambda +i\Omega)\bar z-\Re(z^t\Lambda z)+\cdots
=2y^t\Lambda y+2x^t\Omega y+\cdots,$$ as required.
\end{proof}
\begin{remark}If the discussion of \S\ref{intro} is applied to this normal form
then, at the origin,
$$H=\spn\{\pa_{x_1},\dots,\pa_{x_n},
          \pa_{y_1},\dots,\pa_{y_n}\} \mbox{ and }
N=\spn\{\pa_{y_1},\dots,\pa_{y_n}\}.$$
Upon identifying $H$ with $\C^n$ and $N$ with $\R^n$,
the Levi form of $r$ on $H$
may be identified with the Hermitian matrix $L=\Lambda +i\Omega$,
and its restriction to $N$ with the symmetric matrix $\Lambda$.
\end{remark}
\begin{remark}\label{smoothness}
Save for the initial complex linear change of co\"ordinates,
the change of co\"ordinates constructed in the proof is completely
determined.  It follows easily that locally near a given point of $E$, the
new co\"ordinates, the new defining function for $M$, and the graphing
function for $E$ can be chosen to depend smoothly on $p \in E$.  Moreover,
given $\tau \in J(T_pE)$ transverse to $M$,
by suitably normalizing the initial linear change of co\"ordinates one can
guarantee that $\tau = \pa_v$ in the normal co\"ordinates, with smooth
dependence on $\tau$ as well.
\end{remark}

\section{Proof of Theorem~\ref{main}}
Let $W$ be a two-sided wedge in $M$ as in the statement of
Theorem~\ref{main}.  We will show that $\widehat W$ contains a neighborhood
of $p_0$ in $\C^{n+1}$.  Theorem~\ref{main} follows from this and from an
extension of the Baouendi-Treves approximation theorem \cite{bt} to the
setting of functions
defined on a two-sided wedge in $M$ with a maximally real edge.
This extension states that given a point on the edge of such a wedge
and a direction of the wedge at that point,
there is a two-sided subwedge having the given direction at the given point
and on which any continuous CR function on the wedge, continuous
across its edge, may be uniformly approximated by a sequence of
polynomials.  Since the subwedge also satisfies the hypotheses of
Theorem~\ref{main}, its hull contains a neighborhood of $p_0$, so
the sequence of approximating polynomials also converges on this
neighborhood and therefore defines the desired extension.
The validity of the extended Baouendi-Treves theorem can be seen in either
of two ways.  Rosay \cite{r} has observed that the proof
of the usual Baouendi-Treves theorem applies to
holomorphic functions on wedges in $\C^{n+1}$ continuous across the edge,
and the same observation holds for CR functions, or, more
generally, solutions of the involutive structure on a hypo-analytic
manifold, defined on wedges and continuous across a maximally real edge.
Alternatively, one may deduce the extended theorem
by applying the usual Baouendi-Treves theorem on the
blow-up of $M$ along $E$.  This blow-up has a natural hypo-analytic structure
induced from the embedded CR structure on $M$, and any continuous CR
function on a two-sided wedge in $M$ lifts as a continuous solution of the
underlying involutive structure on the blow-up on a full open set.
(The hypo-analytic structure on the blow-up is discussed in~\cite{micro}.)

Without loss of generality, we may shrink $M$ to suppose that its Levi form is
indefinite throughout~$M$.  Recall then that
by the Lewy extension theorem, any CR function
defined on $M$ uniquely extends as a holomorphic function to a
neighborhood thereof in~$\C^{n+1}$. Also, the usual proof as in \cite{h}
uses an analytic family of discs whose boundaries
are contained in a small neighborhood of a point on $M$.
{From} the existence of these discs one obtains an estimate on the polynomial
hull of subsets
of $M$.  Let $B(z,r)$ denote the usual Euclidean ball of radius $r$
centered at $z \in \C^{n+1}$.  Then from the Lewy discs it follows that
for any relatively compact open subset $U\subset \overline U\subset M$
there are
$\alpha,r_0>0$ such that if $0<r<r_0$ and $z \in U$, then
\begin{equation}\label{balls}
B(z,\alpha r^2) \subset \widehat{(B(z,r) \cap M)}.
\end{equation}

As previously noted, we may assume that the axis $\sigma(p_0)$ of the
defining cone at $p_0$ for the wedge $W$ in $M$ is null.
For $z \in W$, let $r(z)$ denote the distance from $z$ to $E$.  By halving
the aperture and shortening the extent of the cones used in defining $W$, we
obtain another two-sided wedge $\widetilde W$ in $M$ with the same axes as
$W$ and with the property that there is a constant $\kappa>0$ such that
$M \cap B(z,\kappa r(z)) \subset W^{\pm}$ for all $z \in \widetilde
W^{\pm}$.
Combining this with (\ref{balls}) and absorbing $\kappa$ into $\alpha$,
it follows that for some $\alpha >0$,
\begin{equation}\label{Lewy}
\bigcup_{z\in \widetilde W^{\pm}} B(z,\alpha r(z)^2) \subset
\widehat{W^{\pm}}.
\end{equation}
If we could replace $r(z)^2$ by $r(z)$ in (\ref{Lewy}), then we would have
a wedge in $\C^{n+1}$ and we could apply the classical edge-of-the-wedge
theorem to finish the proof.  Roughly speaking, the rest of
the proof consists of using Lemma~\ref{spike_lemma} to go the extra distance.

For $p \in E$ near $p_0$, we can choose normal co\"ordinates as in
Lemma~\ref{normal} and Remark~\ref{smoothness}
depending smoothly on $p$, and we can certainly arrange
that the null axis $\sigma(p_0)$ of the cone defining $W^+$
at $p_0$ points in
the positive $y_1$-direction at the origin in the normal co\"ordinates for
$p_0$.   In one of these normal co\"ordinate systems, we have
\begin{equation}\label{E}
E=\{y=f(x,u),v=g(x,u)\}
\end{equation}
and
\begin{equation}\label{M}
M=\{v=y^t\Lambda y+x^t\Omega y + \phi(x,y,u)\},
\end{equation}
where
\begin{equation}\label{Eerror}
|f(x,u)|, |g(x,u)| = O(|x|^4 + |u|^4)
\end{equation}
and
\begin{equation}\label{Merror}
|\phi(x,y,u)|=O(|x|^3+|y|^3+|u|^3).
\end{equation}
Since $E\subset M$, it follows easily that also
\begin{equation}\label{M+E}
|\phi(x,0,u)|=O(|x|^4+|u|^4).
\end{equation}

{For} fixed small
$$t=(t_2,t_3,\dots,t_n)\in\R^{n-1}$$
(say $|t|\leq 1$ at least),
consider the linear embedding $\Psi_{t}:\C^2\hookrightarrow\C^{n+1}$
given by
$$\Psi_{t}(\zeta_1,\zeta_2)= \zeta_1(1,t_2,\dots,t_n,0)
                         +\zeta_2(0,0,\dots,0,1).$$
Write $(\zeta_1,\zeta_2)=(\xi_1+i\eta_1,\xi_2+i\eta_2)=\xi + i\eta$ and
$\Pi_{t}= \Psi_{t}(\C^2) \subset \C^{n+1}$.
We may compute how $\Pi_{t}$ intersects~$M$.
On $\Pi_{t}$,
$$v-y^t\Lambda y-x^t\Omega y=\eta_2-Q(t)\eta_1{}^2$$
where
$$Q(t)=(1,t)\;\Lambda
\raisebox{-1pt}{$\left(\!\!\begin{array}c
1\\ t\end{array}\!\!\right).$}$$
Thus from (\ref{M}) it follows that $\Pi_{t}\cap M$
is given by the equation
\begin{equation}
\eta_2=\chi(t,\xi,\eta_1)
\end{equation}
where
$\chi(t,\xi,\eta_1) = Q(t)\eta_1{}^2 +
\phi(\xi_1(1,t),\eta_1(1,t), \xi_2).$
Taylor expanding in $\eta_1$, we can write
\begin{equation}\label{chi}
\chi(t,\xi,\eta_1) = Q(t)\eta_1{}^2 + a(t,\xi)+b(t,\xi)\eta_1
+c(t,\xi)\eta_1{}^2 + d(t,\xi,\eta_1)\eta_1{}^3
\end{equation}
for smooth functions $a,b,c,d$, and by (\ref{Merror}) and (\ref{M+E})
we have for some constant $A>0$ and for $k=0,1,2$
\begin{equation}\label{abcd}
|a(t,\xi)|\leq A|\xi|^4,\;|b(t,\xi)|\leq A|\xi|^2,\;
|c(t,\xi)|\leq A|\xi|,\;
|\pa_{\eta_1}^kd(t,\xi,\eta_1)|\leq A.
\end{equation}

For given $(t,\xi)$, we denote by
$(Y(t,\xi),V(t,\xi))\in \R^n \times \R$ the
imaginary part of the unique point of $E$ whose real part is
$\Psi_{t}(\xi_1,\xi_2) = (\xi_1(1,t),\xi_2)$.  {From} (\ref{E}) we have
$Y(t,\xi)=f(\xi_1(1,t),\xi_2)$ and
$V(t,\xi)=g(\xi_1(1,t),\xi_2)$.
{From} (\ref{Eerror}) it follows that there is a constant $B>0$ so that
\begin{equation}\label{Y}
|Y(t,\xi)|\leq B|\xi|^4,\;|V(t,\xi)|\leq B|\xi|^4.
\end{equation}
We sometimes write $Y=(Y_1,Y')$, where $Y'=(Y_2,\ldots,Y_n)$.

The wedge $\widetilde W^+$ in (\ref{Lewy}) has axis
$\sigma(p_0)$ at $p_0$.
As we chose our co\"ordinates to make this the positive $y_1$-axis in
the normal co\"ordinates centered at $p_0$, it follows that there is a fixed
round cone $C$ in $\R^{n+1}$ with axis in the positive $y_1$-direction so
that for all $p$ sufficiently close to $p_0$, $M \cap (E+iC) \subset
\widetilde W^+$ sufficiently near the origin
in the normal co\"ordinates at $p$.
\begin{lemma}\label{Minwedge}
There are constants $\epsilon, \gamma, K >0$ so that if $|t| \leq
\gamma$ and
\begin{equation}\label{partofM}
\begin{array}l K|\xi|^4 < \eta_1 - Y_1(t,\xi)\leq \epsilon,\\[3pt]
\eta_2 = \chi(t,\xi,\eta_1),
\end{array}
\end{equation}
then $\Psi_{t}(\zeta_1,\zeta_2) \in \widetilde W^+$.
\end{lemma}
\begin{proof}
First observe that (\ref{partofM}) and(\ref{Y}) imply that
\begin{equation}\label{xi&eta}
|\xi|^4 \leq \epsilon/K \mbox{ and }|\eta_1|\leq \epsilon(1+B/K),
\end{equation}
so it follows easily upon choosing $\epsilon$ small enough that
$(\zeta_1,\zeta_2)$ lies in a small neighborhood of the origin.
The second line of (\ref{partofM}) implies that
$\Psi_{t}(\zeta_1,\zeta_2) \in M$, so it suffices to show that
$\Psi_{t}(\zeta_1,\zeta_2) \in E+iC$.  If $\delta$ denotes the aperture
of $C$, this follows from the inequalities
$$|\eta_1 t - Y'(t,\eta)|< \delta/2(\eta_1-Y_1(t,\xi))$$
and
$$|\eta_2-V(t,\xi)|< \delta/2(\eta_1-Y_1(t,\xi)).$$
Now
$$|\eta_1 t - Y'(t,\xi)|\leq |(\eta_1-Y_1(t,\xi))t|+
|Y_1(t,\xi)t|+|Y'(t,\xi)|$$
$$\leq \gamma (\eta_1-Y_1(t,\xi)) + 2B|\xi|^4$$
$$< (\gamma+2B/K)(\eta_1-Y_1(t,\xi)),$$
so the first inequality is satisfied if $\gamma\leq \delta/4$ and
$K\geq 8B/\delta$.

For the second inequality, we will need to estimate
$$\chi(t,\xi,\eta_1)=a(t,\xi)+[b(t,\xi) + (Q(t)+c(t,\xi))
\eta_1 +d(t,\xi,\eta_1)\eta_1{}^2]\eta_1.$$
Observe first from (\ref{abcd}) and (\ref{xi&eta}) that
$$|b(t,\xi) + (Q(t)+c(t,\xi))\eta_1 +
d(t,\xi,\eta_1)\eta_1{}^2|\leq \tilde{A}\sqrt{\epsilon},$$
where $\tilde{A}$ is a constant depending only on
$A,B$ and a bound for the matrix $\Lambda $.  Therefore we obtain
$$|\eta_2-V(t,\xi)|\leq |\chi(t,\xi,\eta_1)|+|V(t,\xi)|$$
$$\leq |a(t,\xi)|+\tilde{A}\sqrt{\epsilon}|\eta_1| + B|\xi|^4$$
$$\leq (A+B)|\xi|^4 + \tilde{A}\sqrt{\epsilon}(\eta_1-Y_1(t,\xi))
+\tilde{A}\sqrt{\epsilon}|Y_1(t,\xi)|$$
$$\leq (A+B+\tilde{A}B\sqrt{\epsilon})|\xi|^4 +
\tilde{A}\sqrt{\epsilon}(\eta_1-Y_1(t,\xi))$$
$$<[\tilde{A}\sqrt{\epsilon}
+ (A+B+\tilde{A}B\sqrt{\epsilon})/K](\eta_1-Y_1(t,\xi)).$$
If we choose $\epsilon$ so that $\tilde{A}\sqrt{\epsilon}\leq \delta/4$
and $K$ so that
$$K\geq 4(A+B+\tilde{A}B\sqrt{\epsilon})/\delta,$$
then the second inequality holds as well.
\end{proof}

We now insist that $|t|\leq \gamma$ and we
apply (\ref{Lewy}) to the points $z=\Psi_{t}(\zeta_1,\zeta_2) \in
\widetilde W^+$ of Lemma~\ref{Minwedge}.  Observe that for such
$z$, $r(z)$ is comparable to $\eta_1-Y_1(t,\xi)$.
We deduce that there is $\alpha>0$ so that if we define
$$S_{t}=
\left\{\begin{array}l(\zeta_1,\zeta_2):
K|\xi|^4 < \eta_1 - Y_1(t,\xi)\leq \epsilon \;\mbox{\rm\ and}\\[3pt]
\phantom{(\zeta_1,\zeta_2):{}}
|\eta_2-\chi(t,\xi,\eta_1)|\leq \alpha (\eta_1 - Y_1(t,\xi))^2
\end{array}\right\},$$
then $\Psi_{t}(S_{t}) \subset \widehat{W^+}$.  It follows that
$\Psi_{t}(\widehat{S_{t}}) \subset \widehat{W^+}$, where
$\widehat{S_{t}}$ denotes the hull in $\C^2$.  We intend to apply
Lemma~\ref{spike_lemma} to $S_{t}$.  To this end, define
$$q_1(t,\xi)=Y_1(t,\xi) +  K|\xi|^4,$$
$$q_2(t,\xi)=\chi(t,\xi,q_1(t,\xi)),$$
$$m(t,\xi)=(\pa_{\eta_1}\chi)(t,\xi,q_1(t,\xi)).$$
\begin{lemma}\label{spikesin}
There are positive constants $\beta,\ell,r,\hat{\alpha}$ so that if
$|Q(t)|\leq \hat{\alpha}$, then
$$S_{t} \supset
\left\{\begin{array}l(\zeta_1,\zeta_2):
|\xi|<r\;\mbox{\rm\ and}\\[3pt]
\phantom{(\zeta_1,\zeta_2):{}}(\eta_1,\eta_2)\in
Sp(\beta,\ell,q_1(t,\xi),q_2(t,\xi),m(t,\xi))
\end{array}\right\}.$$
\end{lemma}
\begin{proof}
If $|\xi|<r$ and $(\eta_1,\eta_2)\in
Sp(\beta,\ell,q_1(t,\xi),q_2(t,\xi),m(t,\xi))$,
then $0<\eta_1-q_1(t,\xi)<\ell$, which implies that
$$K|\xi|^4<\eta_1-Y_1(t,\xi)<\ell+K|\xi|^4\leq \ell+Kr^4.$$
Therefore we can guarantee the first condition defining $S_{t}$ by
choosing $\ell\leq \epsilon/2$ and $r$ so small that $Kr^4 \leq \epsilon/2$.
For the second condition, we have
$$|\eta_2 - \chi(t,\xi,\eta_1)|\leq |\eta_2-q_2(t,\xi) -m(t,\xi)
(\eta_1-q_1(t,\xi))|$$
$$+|\chi(t,\xi,\eta_1)-[\chi(t,\xi,q_1(t,\xi))
+(\pa_{\eta_1}\chi)(t,\xi,q_1(t,\xi))(\eta_1-q_1(t,\xi))]|.$$
By definition of $Sp(\beta,\ell,q_1,q_2,m)$,
the first term is at most $\beta(\eta_1-q_1(t,\xi))^2$.
The second term is equal to
$$\left|\int_{q_1(t,\xi)}^{\eta_1} \pa_{\lambda}^2 \chi(t,\xi,\lambda)
(\eta_1-\lambda)d\lambda\right|.$$
Now (\ref{chi}) and (\ref{abcd}) give
$$|\pa_{\lambda}^2 \chi(t,\xi,\lambda)|=|2(Q(t)+c(t,\xi))
+\pa_{\lambda}^2(d(t,\xi,\lambda)\lambda^3)|$$
$$\leq 2|Q(t)|+2A|\xi|+A'|\lambda|$$
for some constant $A'$, so the second term is at most
$$(|Q(t)|+A|\xi|+A'(|\eta_1|+|q_1(t,\xi)|))(\eta_1-q_1(t,\xi))^2.$$
Collecting the terms, we deduce that
$$|\eta_2 - \chi(t,\xi,\eta_1)|\leq [\beta + |Q(t)|
+A|\xi| +A'|\eta_1| +A'(B+K)|\xi|^4](\eta_1-q_1(t,\xi))^2.$$
Since $|\xi|$ and $|\eta_1|$ can be made small by choosing $\ell$ and $r$
small, the result follows.
\end{proof}

Now $Q(0)=0$ in the normal co\"ordinates for $p_0$
since the $y_1$-axis is null in that case.  By continuity it follows that
for all $p$ sufficiently close to $p_0$ and sufficiently small $t$ we
have $|Q(t)|\leq \hat\alpha$.
{From} (\ref{Y}) we obtain $|q_1(t,\xi)| = O(|\xi|^4)$.  Then (\ref{chi})
and (\ref{abcd}) show that $|q_2(t,\xi)| = O(|\xi|^4)$ and
$|m(t,\xi)| = O(|\xi|^2)$.  Therefore we can apply
Lemma~\ref{spike_lemma}
to deduce that
that there is a $\delta>0$ so that
$\Psi_{t}(T_{\delta}) \subset \widehat{W^+}$.  As $t$ varies, the
$\Psi_{t}(T_{\delta})$ sweep out a cone in $i\R^{n+1}$ with axis
in the $y_1$-direction.  As $p$ varies, these cones in the original
variables
sweep out a one-sided wedge in $\C^{n+1}$.  The same argument applies to
$W^-$, giving a two-sided wedge in $\C^{n+1}$ contained in $\widehat W$.
Finally, Rosay~\cite{r} has shown that the hull of such a two-sided wedge in
$\C^{n+1}$ contains a neighborhood of a point on the edge, concluding the
proof of Theorem~\ref{main}.

\begin{remark}\label{sameaxis}
Up to the last step, the construction of the analytic discs in the above
argument proceeds separately on the two sides of the wedge.  Therefore the
argument establishes the following result about extension from a one-sided
wedge.  Let $E\subset M$ be a
maximally real submanifold of a hypersurface $M \subset \C^{n+1}$, suppose
that the Levi form of $M$ at $p_0 \in E$ is indefinite, and let $W$ be a
one-sided wedge in $M$ with edge $E$ having a null direction
$\sigma \in N_{p_0}$.  Then there is a one-sided wedge $\Wp$ in $\C^{n+1}$
near $p_0$ with edge $E$ and axis $\sigma$ at $p_0$ with the property that
every CR function on $W$, continuous in $W\cup E$, extends to a holomorphic
function in $\Wp$, continuous in $\Wp \cup E$.  From this result one can
derive versions of the edge-of-the-wedge theorem in $M$ from those in
$\C^{n+1}$ for extension from unions of one-sided wedges which are not
opposite.
\end{remark}

\section{Proof of Theorem~\ref{one-sided}}
Let $W$ be a two-sided wedge in $M$ with edge $E$ as in the statement of
Theorem~\ref{one-sided}.  We will show that $\widehat{W\cup E}$ contains
an open set $U\subset\{r<0\}$ with properties \ref{firstitem}
and~\ref{seconditem}.  The result then follows upon application of the
Baouendi-Treves approximation theorem as in the proof of Theorem~\ref{main}.

For $p$ near $p_0$ and $\tau \in J(T_pE)$ such that
$dr(\tau)<0$, we may choose normal co\"ordinates as in Lemma~\ref{normal}
and Remark~\ref{smoothness}.  The defining function in the statement of
Theorem~\ref{main} may be replaced by that in Lemma~\ref{normal} since they
are positive smooth multiples of one another.
In the normal co\"ordinates, $\sigma$ is represented by a vector in the
$y$-plane satisfying $\sigma^t \Lambda \sigma >0$.  By scaling we may
suppose that $\sigma^t \Lambda \sigma = 2$.  We may write $E$ and $M$
in the forms (\ref{E}), (\ref{M}) such that (\ref{Eerror}) and
(\ref{Merror}) hold.

We now make the further co\"ordinate change
$$\hat z=z \qquad \hat w=w+\frac{1}{2}iz^t\Lambda z.$$
To first order at the origin this is the identity but $M$ is now defined by
the vanishing of
\begin{equation}\label{Lewyform}
\hat r(\hat x,\hat y,\hat u)=-\hat v+ \frac{1}{2}\hat z^t(\Lambda+i\Omega)
\overline{\hat z}+\hat{\phi}(\hat x,\hat y,\hat u)
\end{equation}
and $E$ by equations of the form
$$
\hat y =\hat f(\hat x,\hat u),\qquad \hat v=\frac{1}{2}\hat x^t \Lambda \hat
x+\hat g(\hat x,\hat u),
$$
where
$$|\hat{\phi}(\hat x,\hat y,\hat u)|=
O(|\hat x|^3 + |\hat y|^3 + |\hat u|^3) \mbox{ and }
|\hat f|,|\hat g|=O(|\hat x|^4 +|\hat u|^4).
$$
The form (\ref{Lewyform}) is that usually adopted in the proof
of Lewy extension~\cite{h}.  The usual Lewy discs are obtained for
small $\delta>0$ by intersecting $\{\hat r < 0\}$ with the image of
the holomorphic embedding
$\C\hookrightarrow\C^{n+1}$ given by
\begin{equation}\label{embed}
\xi+i\eta=\zeta\longmapsto(\hat z,\hat w)=
(\zeta\sigma,i\delta^2).
\end{equation}
This intersection corresponds to the subset $\Delta_{\delta} \subset \C$
given by
$$|\zeta|^2 + \hat{\phi}(\xi\sigma,\eta\sigma,0) < \delta^2,$$
whose connected component containing the origin is a star-shaped domain
with smooth boundary.

Consider first the case $\hat f = \hat g = 0$, which can be arranged if $E$
is real-analytic.  Then the images of $\zeta = \pm \delta$ lie on $E$.
Therefore $\zeta = \pm \delta \in \pa \Delta_{\delta}$.  Since the
$\hat y$-components
of the images of all other points in $\pa \Delta_{\delta}$ are
multiples of $\sigma$, it follows that these boundary points lie in
$W$ for $\delta$ sufficiently small.  Note that the images as $\delta$
varies of the two pure imaginary points of
$\pa \Delta_{\delta}$ form curves in $W^{\pm}$
which are tangent to ($\pm i\sigma,0)$ at the origin.
Since the boundary of the disc is completely contained in $W\cup E$,
the disc itself is contained in
$\widehat{W\cup E}$.  The image of $\zeta = 0$ is the point
$(0,i\delta^2)$, so upon varying $\delta$ we obtain a line segment in the
$v$-direction contained in $\widehat{W\cup E}$.  In the original variables
this corresponds to a curve emanating from
$p$ in the direction $\tau$.  Now let $\tau$, $p$, and $\sigma$ vary.
{From} the above properties it follows easily that the union of the
resulting discs contains an open set $U$ with the required properties.

The above discs must be modified in case $\hat f$ or $\hat g$ do not
vanish, since in that case
they might miss $E$.  In general, the points
$$(\pm\delta\sigma+i\hat f(\pm\delta\sigma,0),i(\delta^2+\hat g
(\pm\delta\sigma,0)))$$ lie on $E$.
Modify the embedding (\ref{embed}) by sending $\zeta=\pm\delta$ to
these two points and uniquely extend to be of the form $\zeta\mapsto\zeta A+B$.
Then it can be shown that for sufficiently small $\delta$
the intersection of $\{\hat r <0\}$ with the image of this embedding
is an analytic disc with boundary in
$W\cup E$ crossing between $W^+$ and $W^-$ at just these two points on~$E$.
The image of $\zeta = 0$ no longer lies on the $v$-axis, but upon varying
$\delta$ we still obtain
a curve tangent to it at the origin, so upon varying $\tau$, $p$, and
$\sigma$ we still obtain a set $U$ as before.

\end{document}